\documentclass[11pt]{article}

\usepackage{amsfonts,amsmath,latexsym,color,epsfig,hyperref}
\newtheorem{theorem}{Theorem}
\newtheorem{lemma}{Lemma}

\def\qed{\ifvmode\mbox{ }\else\unskip\fi\hskip 1em plus 10fill$\Box$}

\input{epsf}
\makeatletter
\def\Ddots{\mathinner{\mkern1mu\raise\p@
\vbox{\kern7\p@\hbox{.}}\mkern2mu
\raise4\p@\hbox{.}\mkern2mu\raise7\p@\hbox{.}\mkern1mu}}
\makeatother

\title{\vspace{-0.7cm}A sequence of triangle-free pseudorandom graphs}
\author{David Conlon\thanks{Mathematical Institute, Oxford OX2 6GG,
United Kingdom. Email: {\tt david.conlon@maths.ox.ac.uk}. Research
supported by a Royal Society University Research Fellowship and by ERC Starting Grant 676632.}}
\date{}

\begin{document}
\maketitle

\begin{abstract}
A construction of Alon yields a sequence of highly pseudorandom triangle-free graphs with edge density significantly higher than one might expect from comparison with random graphs. We give an alternative construction for such graphs.
\end{abstract}

\section{Introduction}

A graph $G$ is said to be $(p, \beta)$-jumbled if
\[\left|e(X) - p \binom{|X|}{2}\right| \leq \beta |X|\]
for all $X \subseteq V(G)$. For example, the binomial random graph $G_{n,p}$ is $(p, \beta)$-jumbled with $\beta = O(\sqrt{pn})$. It is not hard to show~\cite{EGPS88, ES71} that this is essentially best possible, in that a graph with $n$ vertices cannot be $(p, \beta)$-jumbled with $\beta = c\sqrt{pn}$ for $c$ sufficiently small. For further information on jumbled graphs and their properties, we refer the reader to the survey~\cite{KS06} or, for more recent developments, the paper~\cite{CFZ14}.

One important class of $(p, \beta)$-jumbled graphs is the collection of $(n, d, \lambda)$-graphs. These are $d$-regular graphs on $n$ vertices such that all eigenvalues of the adjacency matrix, save the largest, are bounded in absolute value by $\lambda$. By the famous expander mixing lemma, these graphs are $(p, \beta)$-jumbled with $p = d/n$ and $\beta = \lambda$. 

One of the best known examples of a pseudorandom graph, constructed by Alon~\cite{A94}, is a triangle-free $(n, d, \lambda)$-graph with $n = 2^{3k}$, $d = 2^{k-1} (2^{k-1} - 1)$  and $\lambda = O(2^k)$. Taking $p = d/n$, we have $\sqrt{p n} = \sqrt{d} = \Omega(2^k)$, so the graph is close to optimally pseudorandom. Since $p = \Omega(n^{-1/3})$, the construction also has surprisingly high density. While there are various ways to modify the usual random graph to produce triangle-free graphs with density roughly $n^{-1/2}$ (see, for example,~\cite{B09}), no such modification can hope to push very far past this density. Nevertheless, Alon's construction does so. The purpose of this note is to give another construction for such graphs.

\begin{theorem} \label{thm:main}
There exists a sequence $(n_i)_{i=1}^{\infty}$ of positive integers such that, for each $i \geq 1$, there is a triangle-free graph $G_i$ on $n_i$ vertices which is $(p, \beta)$-jumbled with $p = \Omega(n_i^{-1/3})$ and $\beta = O(\sqrt{p n_i} \log n_i)$.
\end{theorem}

Our construction is weaker than Alon's on several counts: it does not produce regular graphs; it is not completely explicit; and it does not generalise easily. One might also level the accusation that the resulting graphs are not optimally pseudorandom, with the condition $\beta = O(\sqrt{pn} \log n)$ being a logarithmic factor away from the desired bound. However, it seems likely that this extra log factor is simply an artifact of our proof. Countering these disadvantages, we believe that our construction, which we describe in detail below, is more intuitive than Alon's.

For concreteness, we will work with the polarity graph of Lazebnik, Ustimenko and Woldar~\cite{LUW99}, though the role played by this graph could also be taken by a number of other $C_6$-free graphs with $n$ vertices and $\Omega(n^{4/3})$ edges. Suppose then that $q$ is an odd power of $2$ and let $n = q^3 + q^2 + q + 1$. The polarity graph, which has a small number of loops, is an $(n, d, \lambda)$-graph with $d = q + 1 \geq n^{1/3}$ and, as noted in~\cite[Section 3.7]{KS06}, $\lambda = \sqrt{2q} = O(n^{1/6})$. Once the loops are removed, the resulting graph, which we label by $H$, is $C_3, C_4$ and $C_6$ free. For each vertex $v$ in $H$, we randomly partition its neighbourhood $N_H(v)$ into two sets $A_v$ and $B_v$ and let $G_v$ be the complete bipartite graph between $A_v$ and $B_v$. We now define $G$ to be the graph with the same vertex set as $H$ and edge set $\cup_{v \in V(H)} G_v$. We will show that $G$ asymptotically almost surely satisfies the requirements of Theorem~\ref{thm:main}.

It is straightforward to verify that $G$ contains no triangles. To begin, note that since $H$ is $C_4$-free, the only edges of $G$ in $N_H(v)$ are those in $G_v$. Now suppose that $abc$ is a triangle in $G$. If $a, b$ and $c$ are all contained in the same neighbourhood $N_H(v)$ they cannot form a triangle, since the only edges of $G$ in $N_H(v)$ are those in $G_v$ and $G_v$ is bipartite. It must then be the case that the edges $ab$, $bc$ and $ca$ are contained in three different neighbourhoods $N_H(u)$, $N_H(v)$ and $N_H(w)$, respectively. If $u, v$ and $w$ are all distinct from $a, b$ and $c$, then $a u b v c w$ would form a cycle of length $6$ in $H$, contradicting the fact that $H$ is $C_6$-free. On the other hand, if $u = c$, say, then $H$ must contain the triangle $b v c$, again contradicting the choice of $H$. 

The remaining claim, that $G$ is asymptotically almost surely $(p, \beta)$-jumbled with $p = \Omega(n^{-1/3})$ and $\beta = O(\sqrt{p n} \log n)$, will be verified in the next section.

\section{Proving $G$ is jumbled} \label{sec:jumble}

Let $H_0$ be the polarity graph with $n = q^3 + q^2 + q + 1$ vertices. This graph is $(q+1)$-regular, has $q^2 + 1$ loops and all eigenvalues of the adjacency matrix, save the largest, are bounded in absolute value by $\sqrt{2q}$. We will form $G$ from $H_0$ by a slightly different procedure to that described in the introduction, though the two are easily seen to produce the same graph.

In the first step, we form $H_0^{(2)}$, the multigraph with loops on the same vertex set as $H_0$ where two vertices are joined if there is a walk of length two in $H_0$ between them, allowing for multiple edges if there is more than one such walk. We note that each vertex has $q+1$ loops, one for each edge in $H_0$, and, since $H$, the simple graph formed by removing the loops from $H_0$, is $C_3$ and $C_4$-free, the only parallel edges arise from loops in $H_0$. In the next step, we turn $H_0^{(2)}$ into a simple graph by removing all loops from $H_0^{(2)}$ and all edges whose corresponding walk in $H_0$ used a loop. The resulting graph $G_1$ is easily seen to be the union of $n$ cliques, each clique being $N_H(v)$ for some $v$ in $H$. We now form the required graph, as before, by randomly partitioning $N_H(v)$ into two sets $A_v$ and $B_v$ and letting $G$ be the union over all $v$ in $H$ of the complete bipartite graphs between $A_v$ and $B_v$.

Following this plan, we first look at $H_0^{(2)}$. Letting $M$ be the adjacency matrix of $H_0$, the adjacency matrix of $H_0^{(2)}$ is simply $M^2$, which implies that the eigenvalues of $H_0^{(2)}$ are the squares of the eigenvalues of $H_0$. Therefore, $H_0^{(2)}$ is an $(n, d, \lambda)$-graph with $d = (q + 1)^2$ and $\lambda = 2q$. Since each vertex of $H_0^{(2)}$ is contained in exactly $q + 1$ loops, the graph $G_0$ formed by removing these loops has adjacency matrix $M^2 - (q+1)I$, implying that $G_0$ is an $(n, d, \lambda)$-graph with $d = q(q+1)$ and $\lambda = q + 1$. Therefore, by the expander mixing lemma, for all $X \subseteq V(G_0)$, 
\[\left|e_{G_0}(X) - \frac{q}{q^2+1} \binom{|X|}{2}\right| \leq (q + 1) |X|,\]
where we used the fact that $n = q^3 + q^2 + q + 1 = (q^2+1)(q+1)$.

Note now that every loop in $H_0$ has exactly $q$ other neighbours with which it can form a non-degenerate walk of length two. Therefore, since any $X \subseteq V(H_0)$ contains at most $|X|$ loops, we remove at most $q |X|$ edges from $X$ when forming $G_1$ from $G_0$. By the estimate above, this implies that, for all $X \subseteq V(G_1)$,
\[\left|e_{G_1}(X) - \frac{q}{q^2 + 1}\binom{|X|}{2}\right| \leq (2q + 1) |X|.\]
Therefore, $G_1$ is $(p, \beta)$-jumbled with $p = q/(q^2 + 1)$ and $\beta = 2q + 1$.

Recall that $G_1$ is the union of cliques, while $G$ is the union of random bipartite graphs, one for each clique in $G_1$. Our aim now is to show that asymptotically almost surely $G$ is $(p, \beta)$-jumbled with $p = q/2(q^2 + 1) = \Omega(n^{-1/3})$ and $\beta = O(q \log n) = O(\sqrt{pn} \log n)$. We will do this in a rather naive fashion, estimating the probability that
\begin{equation} \label{eqn:main}
\left|e_{G}(X) - \frac{q}{2(q^2 + 1)}\binom{|X|}{2}\right| \leq C q |X| \log n
\end{equation}
for any given $X$ and taking a union bound. To do this, we will need the following concentration inequality for quadratic forms in independent random variables due to Hanson and Wright~\cite{HW71} (the exact version we state follows from Theorem 1.1 in~\cite{RV13}).

\begin{lemma} \label{lem:HW}
Let $Z = (Z_1, \dots, Z_t) \in \{-1, +1\}^t$ be a random vector with independent components each of which is equal to $1$ or $-1$ with probability $1/2$. Let $M$ be a $t \times t$ real matrix. Then
\[\mathbb{P}[|Z^T M Z - \mathbb{E}(Z^T M Z)| > \epsilon] \leq 2 \exp\left\{- c \min \left(\frac{\epsilon^2}{\|M\|_{F}^2}, \frac{\epsilon}{\|M\|}\right)\right\}, \]
where $\|M\|_{F} = (\sum_{i, j} m_{ij}^2)^{1/2}$ is the Frobenius norm and $\|M\| = \sup_{x \neq 0} \|M x\|_2/\|x\|_2$ is the spectral norm.  
\end{lemma}

Suppose now that $G[X]$ is the union of $s$ cliques $T_1, \dots, T_s$, of orders $t_1, \dots, t_s$, and let $t = t_1 + \dots + t_s$. We define $t$ random variables $X_1, \dots, X_t$, each equal to $1$ or $-1$ with probability $1/2$, and assign one of these random variables to every vertex of every clique, noting that any given vertex may receive multiple random variables, but only one relative to any given clique. 

Suppose that the random variables assigned to the clique $T_i$ are $Z_{i1}, \dots, Z_{it_i}$. If $v_i$ is the vertex whose neighbourhood in $H$ is $T_i$, the value of $Z_{ij}$ determines whether its corresponding vertex $T_i(j)$ is placed in $A_{v_i}$ or $B_{v_i}$, with $T_i(j)$ placed in $A_{v_i}$ if $Z_{ij} = 1$ and $B_{v_i}$ if $Z_{ij} = -1$. The number of edges in $G[T_i]$ is then
\[e_G(T_i) = |A_{v_i}||B_{v_i}| = \left(\frac{t_i}{2} + \frac{1}{2}\sum_{j=1}^{t_i} Z_{ij}\right) \left(\frac{t_i}{2} - \frac{1}{2}\sum_{j=1}^{t_i} Z_{ij}\right) = \frac{t_i^2}{4} - \frac{1}{4} \sum_{j=1}^{t_i} \sum_{k=1}^{t_i} Z_{ij} Z_{ik}.\]
Summing over $i$, we have
\begin{align*}
e_G(X) = \sum_{i=1}^s e_G(T_i) & = \sum_{i=1}^s \frac{t_i^2}{4} - \frac{1}{4} \sum_{i=1}^s \sum_{j=1}^{t_i} \sum_{k=1}^{t_i} Z_{ij} Z_{ik}\\
& = \frac{1}{2} \sum_{i=1}^s \binom{t_i}{2} + \frac{1}{4} \sum_{i=1}^s t_i - \frac{1}{4} \sum_{i=1}^s \sum_{j=1}^{t_i} \sum_{k=1}^{t_i} Z_{ij} Z_{ik}\\
& = \frac{e_{G_1}(X)}{2} - \frac{1}{4} \sum_{i=1}^s \sum_{1 \leq j \neq k \leq t_i} Z_{ij} Z_{ik}. 
\end{align*}
Therefore, by our estimate on $e_{G_1}(X)$, it only remains to show that 
\[Q = \sum_{i=1}^s \sum_{1 \leq j \neq k \leq t_i} Z_{ij} Z_{ik}\]
is smaller than $C q |X| \log n$ with sufficiently high probability.

Let $M$ be the $t \times t$ matrix whose entry $m_{jk}$ is equal to $1$ if $j \neq k$ and $j$ and $k$ (each of which represents a vertex associated to a particular clique) are from the same clique $T_i$. Otherwise, we take $m_{jk}$ to be $0$. Then $Q = Z^T M Z$ and it follows from the Hanson--Wright bound that
\[\mathbb{P}[|Q| > C q |X| \log n] \leq 2 \exp\left\{- c \min \left(\frac{C^2 q^2 |X|^2 \log^2 n}{\|M\|_{F}^2}, \frac{C q |X| \log n}{\|M\|}\right)\right\}.\]
But it is straightforward to verify that 
\[\|M\|_F^2 = \sum_{i=1}^s t_i(t_i-1) = 2 e_{G_1}(X) \leq \frac{|X|^2}{q} + 6q|X| \leq 2 \max\left\{\frac{|X|^2}{q}, 6q|X|\right\}\]
and, writing $\rho(A)$ for the spectral radius of a matrix $A$,
\[\|M\| = \sqrt{\rho(M^* M)} = \rho(M) \leq \sup_{x \neq 0} \frac{\|M x\|_{\infty}}{\|x\|_{\infty}} = \max_{1 \leq j \leq t} \sum_{k=1}^t |m_{jk}| \leq q.\]
Therefore, 
\begin{align*}
\mathbb{P}[|Q| > C q |X| \log n] & \leq 2 \exp\left\{- c \min \left(\frac{1}{2}C^2 q^3 \log^2 n, \frac{1}{12}C^2 q |X| \log^2 n, C |X| \log n\right)\right\}\\
& \leq 2 \exp\{- 2 |X| \log n\}
\end{align*}
for $C$ sufficiently large in terms of $c$. Applying the union bound, we see that the probability there exists a set $X$ such that \eqref{eqn:main} fails is at most
\[\sum_{|X| = 1}^n \binom{n}{|X|} 2 e^{- 2 |X| \log n} \leq 2 \sum_{|X| = 1}^n n^{|X|} e^{-2 |X| \log n} \leq 2 \sum_{|X| = 1}^n e^{-|X| \log n} \leq \frac{2}{n-1}.\]
The result follows.

\section{Concluding remarks}

Our construction also extends to give a sequence $(n_i)_{i=1}^{\infty}$ of positive integers such that, for each $i \geq 1$, there is a $C_5$-free graph $G_i$ on $n_i$ vertices which is $(p, \beta)$-jumbled with $p = \Omega(n_i^{-3/5})$ and $\beta = O(\sqrt{p n_i} \log n_i)$. The construction starts with the $C_{10}$-free polarity graph of Lazebnik, Ustimenko and Woldar~\cite{LUW99}, but follows the proof of Theorem~\ref{thm:main} in all other respects. Because we lack optimal constructions of $C_{2\ell}$-free graphs for $\ell \geq 6$, our method does not extend further to give constructions of pseudorandom $C_{2k+1}$-free graphs for $k \geq 3$. As mentioned in the introduction, this is a distinct weakness of our method when compared to Alon's, which does extend to longer odd cycles~\cite{AK98, KS06}.

An alternative method for proving the jumbledness of our construction $G$ might be to estimate the eigenvalues of its adjacency matrix $M$ by using the fact that Tr$(M^{2k})$ is both the sum of the $2k^{th}$ powers of its eigenvalues and the number of walks of length $2k$ in $G$. However, $G$ is constructed by starting from a graph $G_1$ which is a union of cliques and then taking a random bipartite graph within each clique. This process causes an imbalance between odd and even cycles within each clique, deleting all odd cycles, while doubling the proportion of even cycles relative to the density. This makes it difficult to count the number of degenerate walks of length $2k$ without having close control over the counts of degenerate walks of different types in the base graph $G_1$. Nevertheless, it is plausible that this could be done, and may even allow one to save the lost logarithmic factor in $\beta$.

Another definite weakness of our method is that it is not explicit and so, unlike Alon's example, cannot be used to give a constructive lower bound for the off-diagonal Ramsey number $r(3, t)$. It remains to decide whether there is some more explicit method for choosing large bipartite subgraphs of the cliques in $G_1$ which also produces a highly pseudorandom subgraph.

\vspace{3mm}
{\bf Acknowledgements.} I would like to thank the anonymous referee for a number of helpful remarks.

\end{document}